\documentclass[a4paper,11pt]{amsart}

\usepackage{amssymb,mathrsfs,cite,url}

\newtheorem{theorem}{Theorem}[section]
\newtheorem{lemma}[theorem]{Lemma}
\newtheorem{proposition}[theorem]{Proposition}

\theoremstyle{definition}
\newtheorem{remark}[theorem]{Remark}

\numberwithin{equation}{section}

\makeatletter
\@namedef{subjclassname@2010}{%
	\textup{2010} Mathematics Subject Classification}
\makeatother

\DeclareMathOperator{\RE}{Re}
\DeclareMathOperator{\IM}{Im}

\begin{document}

\title[Mean values of logarithmic derivatives of $L$-functions]
{On certain mean values of logarithmic derivatives of $L$-functions and the related density functions}

\author[M. Mine]{Masahiro Mine}
\address{Department of Mathematics\\ Tokyo Institute of Technology\\ 2-12-1 Ookayama, Meguro-ku, Tokyo 152-8551, Japan}
\email{mine.m.aa@m.titech.ac.jp}

\date{}

\begin{abstract}
We study some ``density function'' related to the value-distribution of $L$-functions. 
The first example of such a density function was given by Bohr and Jessen in 1930s for the Riemann zeta-function. 
In this paper, we construct the density function in a wide class of $L$-functions. 
We prove that certain mean values of $L$-functions in the class are represented as integrals involving the related density functions. 
\end{abstract}



\maketitle

\section{Introduction}\label{sec1}
We begin with recalling a classical result on the value-distribution of the Riemann zeta-function $\zeta(s)$ obtained by Bohr and Jessen. 
For any $\sigma>1/2$, let 
\begin{gather*}
G
=\{ s=\sigma+it \mid \sigma>1/2 \} 
\setminus \bigcup_{\rho=\beta+i\gamma} \{ s=\sigma+i\gamma \mid 1/2<\sigma \leq \beta \}, 
\end{gather*}
where $\rho$ runs through all zeros of $\zeta(s)$ with $\beta>1/2$. 
Then we define $\log{\zeta}(s)$ for $s \in G$ by analytic continuation along the horizontal line. 
Fix a rectangle $R$ in the complex plane whose edges are parallel to the coordinate axes, and denote by $\mathcal{V}_\sigma(T,R)$ the Lebesgue measure of the set 
\begin{gather*}
\{ t \in [-T,T] \mid \sigma+it \in G,~ \log{\zeta}(\sigma+it) \in R\}. 
\end{gather*}
Bohr and Jessen \cite{BohrJessen1930, BohrJessen1932} proved that there exists the limit value 
\begin{gather}\label{eqBJ1}
\mathcal{W}_\sigma(R) 
=\lim_{T \to\infty} \frac{1}{2T} \mathcal{V}_\sigma(T,R) 
\end{gather}
for any fixed $\sigma>1/2$. 
They also showed that there exists a non-negative real valued continuous function $\mathcal{M}_\sigma(z)$ such that the formula 
\begin{gather}\label{eqBJ2}
\mathcal{W}_\sigma(R) 
=\int_{R} \mathcal{M}_\sigma(z) \,|dz| 
\end{gather}
holds with $|dz|=(2\pi)^{-1}dxdy$. 
Their study was developed in various ways, for example, Jessen--Wintner \cite{JessenWintner1935}, Borchsenius--Jessen \cite{BorchseniusJessen1948}, Laurin\v{c}ikas \cite{Laurincikas1996a}, and Matsumoto \cite{Matsumoto1989}. 

Matsumoto \cite{Matsumoto1990} generalized limit formula \eqref{eqBJ1} in a quite wide class of zeta-functions, which is now called the \emph{Matsumoto zeta-functions}. 
On the other hand, an analogue of integral formula \eqref{eqBJ2} was obtained only in some restricted cases, for example, the case of Dedekind zeta-functions of finite Galois extensions of $\mathbb{Q}$ \cite{Matsumoto1992}, and automorphic $L$-functions of normalized holomorphic Hecke-eigen cusp forms of level $N$ \cite{MatsumotoUmegaki2019}. 
Thus it is worth studying ``density functions'' such as $\mathcal{M}_\sigma(z)$ for more general zeta- or $L$-functions. 

Kershner and Wintner \cite{KershnerWintner1937} proved analogues of formulas \eqref{eqBJ1} and \eqref{eqBJ2} for $(\zeta'/\zeta)(s)$. 
In this paper, we construct the density functions $M_\sigma(z;F)$ for functions $F(s)$ in a subclass of the Matsumoto zeta-functions and generalize Kershner--Wintner's result.

\section{$L$-functions and the related density functions}\label{sec2}

\subsection{Class of $L$-functions}\label{sec2.1}
We introduce the class $\mathcal{S}_{\mathrm{I}}$ as the set of all functions $F(s)$ represented as Dirichlet series 
\begin{gather*}
F(s) 
=\sum_{n=1}^{\infty} \frac{a_F(n)}{n^s}
\end{gather*}
in some half plane that satisfy the following axioms: 
\begin{enumerate}
\item \emph{Ramanujan hypothesis.}
Dirichlet coefficients $a_F(n)$ satisfy $a_F(n) \ll_\epsilon n^\epsilon$ for every $\epsilon>0$. 
\item \emph{Analytic continuation.}
There exists a non-negative integer $m$ such that $(s-1)^m F(s)$ is an entire function of finite order. 
\item \emph{Functional equation.}
$F(s)$ satisfies a functional equation of the form 
\begin{gather*}
\Lambda_F(s) 
=\omega \overline{\Lambda_F(1-\overline{s})}, 
\end{gather*}
where 
\begin{gather*}
\Lambda_F(s) 
=F(s) Q^s \prod_{j=1}^{r} \Gamma(\lambda_j s+\mu_j), 
\end{gather*}
with some $|\omega|=1$, $Q>0$, $\lambda_j>0$, $\RE(\mu_j) \geq 0$. 
\item \emph{Polynomial Euler product.}
For $\sigma>1$, $F(s)$ is expressed as the infinite product 
\begin{gather*}
F(s)= 
\prod_{p} \prod_{j=1}^{g} \left(1-\frac{\alpha_j(p)}{p^s}\right)^{-1}, 
\end{gather*}
where $g$ is a positive constant and $\alpha_j(p) \in \mathbb{C}$. 
\item \emph{Prime mean square.}
There exists a positive constant $\kappa$ such that 
\begin{gather*}
\lim_{x \to\infty} \frac{1}{\pi(x)} \sum_{p \leq x}|a_F(p)|^2
=\kappa, 
\end{gather*}
where $\pi(x)$ stands for the number of prime numbers less than or equal to $x$. 
\end{enumerate}
The above axioms come from two classes of $L$-functions introduced by Selberg  \cite{Selberg1992} and Steuding \cite{Steuding2007}. 
We see that the class $\mathcal{S}_{\mathrm{I}}$ is just equal to the intersection of these classes, and it is also a subclass of the Matsumoto zeta-functions, see Section 2 of \cite{Steuding2007}. 

Let $N_F(\sigma,T)$ be the number of zeros $\rho=\beta+i\gamma$ of $F(s)$ with $\beta>\sigma$ and $0<\gamma<T$. 
Then for the function $F(s)$ satisfying axioms $\mathrm{(1)}$--$\mathrm{(4)}$, there exists a positive constant $b$ such that for any $\epsilon>0$, 
\begin{gather}\label{eq:WZD}
N_F(T,\sigma) 
\ll_\epsilon T^{b(1-\sigma)+\epsilon} 
\end{gather}
as $T \to\infty$, uniformly for $\sigma \geq1/2$ \cite[Lemma 3]{KaczorowskiPerelli2003}. 
From the proof of \cite{KaczorowskiPerelli2003}, estimate \eqref{eq:WZD} generally holds with $b=4(d_F+3)$, where $d_F$ is the degree of $F$ defined by 
\begin{gather*}
d_F
=2 \sum_{j=1}^{r} \lambda_j. 
\end{gather*}
The constant $b$ is taken smaller in some special cases, for example, Heath-Brown \cite{HeathBrown1977} showed that the Dedekind zeta-functions attached to algebraic number fields of degree $d \geq 3$ satisfy \eqref{eq:WZD} with $b=d$, and Perelli \cite{Perelli1982} obtained it with $b=d_F$ in a subclass of the Selberg class. 

Next, we define the subclass $\mathcal{S}_{\mathrm{II}}$ as the set of all $F(s)$ satisfying axioms $\mathrm{(1)}$--$\mathrm{(5)}$ and the following $\mathrm{(6)}$: 
\begin{enumerate}
\setcounter{enumi}{5}
\item \emph{Zero density estimate.}
There exist positive constants $c$ and $A$ such that 
\begin{gather}\label{eq:ZD}
N_F(T,\sigma)
\ll T^{1-c(\sigma-\frac{1}{2})} (\log T)^A 
\end{gather}
as $T \to\infty$, uniformly for $\sigma \geq 1/2$. 
\end{enumerate}
There are many zeta- or $L$-functions that belong to the class $\mathcal{S}_{\mathrm{I}}$, for instance, the Riemann zeta-functions $\zeta(s)$, Dirichlet $L$-functions $L(s,\chi)$ of primitive characters $\chi$, Dedekind zeta-functions $\zeta_K(s)$, automorphic $L$-functions $L(s,f)$ of normalized holomorphic Hecke-eigen cusp forms $f$ with respect to ${SL}_2(\mathbb{Z})$. 
Furthermore, estimate \eqref{eq:ZD} is proved for $\zeta(s)$ by Selberg \cite{Selberg1946}, for $L(s,\chi)$ by Fujii \cite{Fujii1974}, and for $L(s,f)$ by Luo \cite{Luo1995}, and hence they belong to the subclass $\mathcal{S}_{\mathrm{II}}$.

\subsection{Statements of results}\label{sec2.2}
For an integrable function $f(z)$, we denote its Fourier transform and Fourier inverse transform by 
\begin{gather*}
\widehat{f}(z)= 
f^{\wedge}(z)= 
\int_{\mathbb{C}} f(w) \psi_z(w) \,|dw| 
\quad\text{and}\quad
f^{\vee}(z)= 
\int_{\mathbb{C}} f(w) \psi_{-z}(w) \,|dw|, 
\end{gather*}
respectively, where $\psi_w(z)=\exp(i \RE(z \overline{w}) )$ is an additive character of $\mathbb{C}$ and $|dw|$ is the measure $(2\pi)^{-1}dudv$ for $w=u+iv$. 
According to \cite[Section 9]{IharaMatsumoto2011a} or \cite[Section 5]{IharaMatsumoto2011b}, we then define the class $\Lambda$ as 
\begin{gather*}
\Lambda= 
\{ f \in L^1 \mid \text{$f, \widehat{f} \in L^1 \cap L^\infty$ and $(f^{\wedge})^{\vee}=f$ holds} \}. 
\end{gather*}
We see that any Schwartz function belongs to the class $\Lambda$, and especially, any compactly supported $C^\infty$-function does. 

The first main result of this paper is related to the mean values of $L$-functions. 

\begin{theorem}\label{thmM1}
Let $F\in\mathcal{S}_{\mathrm{I}}$. 
Let $\sigma_1$ be a large fixed positive real number. 
Let $\theta, \delta>0$ be real numbers with $\delta+3\theta<1/2$. 
Let $\epsilon>0$ be a small fixed real number. 
Let $\Phi \in \Lambda$. 
Then there exists a constant $T_{\mathrm{I}}=T_{\mathrm{I}}(F,\sigma_1,\theta,\delta,\epsilon)>0$ such that the following formula 
\begin{gather}\label{eqM1}
\frac{1}{T} \int_{0}^{T} \Phi \left(\frac{F'}{F}(\sigma+it)\right) \,dt 
=\int_{\mathbb{C}} \Phi(z) M_\sigma(z;F) \,|dz|
+E 
\end{gather}
holds for all $T \geq T_{\mathrm{I}}$ and for all $\sigma \in [1-b^{-1}+\epsilon, \sigma_1]$, where $M_\sigma(z;F)$ is a non-negative real valued continuous function uniquely determined from $F(s)$, and the constant $b$ is that in \eqref{eq:WZD}. 
The error term $E$ is estimated as 
\begin{gather}\label{eqM1E}
E
\ll \exp\left(-\frac{1}{4}(\log T)^{\frac{2}{3}\theta}\right) \int_{\Omega} |\widehat{\Phi}(z)| \,|dz| 
+\int_{\mathbb{C} \setminus \Omega} |\widehat{\Phi}(z)| \,|dz|, 
\end{gather}
where the implied constant depends only on $F, \sigma_1, \epsilon$, and 
\begin{gather*}
\Omega= 
\{ z=x+iy \in \mathbb{C} \mid -(\log T)^\delta \leq x, y \leq (\log T)^\delta \}. 
\end{gather*}
Moreover, if $F \in \mathcal{S}_{\mathrm{II}}$, then there exists a constant $T_{\mathrm{II}}=T_{\mathrm{II}}(F,\sigma_1,\theta,\delta)>0$ such that \eqref{eqM1} and \eqref{eqM1E} hold together with $T\geq T_{\mathrm{II}}$ and $\sigma \in [1/2+(\log T)^{-\theta}, \sigma_1]$, where the implied constant depends only on $F$ and $\sigma_1$. 
\end{theorem}

Then, let again $R$ be a rectangle in the complex plane whose edges are parallel to the axes, and define $V_\sigma(T,R;F)$ as the Lebesgue measure of the set of all $t \in [0,T]$ for which $(F'/F)(\sigma+it)$ belongs to $R$. 
Denote by $\nu_k$ the usual $k$-dimensional Lebesgue measure. 
The second result is an analogue of Bohr--Jessen's limit theorem for $(F'/F)(s)$. 

\begin{theorem}\label{thmM2}
Let $F\in\mathcal{S}_{\mathrm{I}}$. 
Let $\sigma$ be fixed with $\sigma>1-b^{-1}$, where the constant $b$ is that in \eqref{eq:WZD}. 
Let $\epsilon>0$ be an arbitrarily small real number. 
Then we have 
\begin{gather}\label{eqM2}
\frac{1}{T} V_\sigma(T,R;F) 
=\int_{R} M_\sigma(z;F) \,|dz| 
+O\left( (\nu_2(R)+1)(\log T)^{-\frac{1}{2}+\epsilon} \right) 
\end{gather}
as $T \to\infty$, where the implied constant depends only on $F, \sigma$, and $\epsilon$. 
Moreover, if $F \in \mathcal{S}_{\mathrm{II}}$, then \eqref{eqM2} holds with any fixed $\sigma>1/2$. 
\end{theorem}

\subsection{Remarks on the related works}\label{sec2.3}
The Riemann zeta-function $\zeta(s)$ is a typical example of the member of the subclass $\mathcal{S}_{\mathrm{II}}$. 
In this case, Theorem \ref{thmM1} is essentially Theorem 1.1.1 of \cite{Guo1996a}, and the density function $M_\sigma(z;\zeta)$ was used to study of the distribution of zeros of $\zeta'(s)$ in \cite{Guo1996b}. 

Theorem \ref{thmM2} is related to the study on the \emph{discrepancy estimates} for zeta-functions. 
Let 
\begin{gather*}
\mathcal{D}_\sigma(T,R)
=\frac{1}{2T} \mathcal{V}_\sigma(T,R) 
-\mathcal{W}_\sigma(R). 
\end{gather*}
We know that $\mathcal{D}_\sigma(T, R)=o(1)$ as $T \to\infty$ by \eqref{eqBJ1}. 
Matsumoto \cite{Matsumoto1987} gave a better upper bound for $\mathcal{D}_\sigma(T,R)$, which was improved by Harman and Matsumoto \cite{HarmanMatsumoto1994}. 
They proved 
\begin{gather*}
\mathcal{D}_\sigma(T,R)
\ll (\nu_2(R)+1) (\log T)^{-A(\sigma)+\epsilon} 
\end{gather*}
for an arbitrarily small $\epsilon>0$, where 
\begin{gather*}
A(x)
= 
\begin{cases}
(x-1)/(3+2x) 
& \text{if  $x>1$}, \\ 
(4x-2)/(21+8x) 
&\text{if $1/2<x \leq 1$}. 
\end{cases}
\end{gather*}
Matsumoto \cite{Matsumoto2007} also generalized this result for Dedekind zeta-functions even in the case of non-Galois extensions. 
We note that $A(x) \leq 1/2$ for any $x>1/2$. 
Though the difference of logarithms and logarithmic derivatives exists, Theorem \ref{thmM2} gives a better estimate on the discrepancy for $(F'/F)(s)$. 

Recently, Ihara and Matsumoto studied density functions such as $\mathcal{M}_\sigma(z)$ more precisely, and named them ``$M$-functions'' for $L$-functions, see \cite{Ihara2008, IharaMatsumoto2011a, IharaMatsumoto2011b, IharaMatsumoto2014}.

\section{Proof of Theorem \ref{thmM1}}\label{sec3}
We begin with considering the case of $\Phi=\psi_z$ in Theorem \ref{thmM1}. 
The following proposition is a key for the proof of the theorem: 

\begin{proposition}\label{propM1}
Let $F(s)$ be a function satisfying axioms $\mathrm{(1)}$--$\mathrm{(4)}$. 
Let $\sigma_1$ be a large fixed positive real number. 
Let $\theta, \delta>0$ be real numbers with $\delta+3\theta<1/2$. 
Let $\epsilon>0$ be a small fixed real number. 
Then there exists a constant $T_{\mathrm{I}}=T_{\mathrm{I}}(F,\sigma_1,\theta,\delta,\epsilon)>0$ such that we have 
\begin{gather}\label{eqpropM1}
\frac{1}{T} \int_{0}^{T} \psi_z\left(\frac{F'}{F}(\sigma+it)\right) \,dt 
=\widetilde{M}_\sigma(z;F) 
+O\left( \exp\left(-\frac{1}{4}(\log T)^{\frac{2}{3}\theta}\right) \right) 
\end{gather}
for all $T \geq T_{\mathrm{I}}$, for all $\sigma \in [1-b^{-1}+\epsilon, \sigma_1]$, and for all $z \in \Omega$, where $\widetilde{M}_\sigma(z;F)$ is a function uniquely determined from $F(s)$. 
The implied constant depends only on $F, \sigma_1$ and $\epsilon$. 
If $F(s)$ further satisfies axiom $\mathrm{(6)}$, there exists a constant $T_{\mathrm{II}}=T_{\mathrm{II}}(F,\sigma_1,\theta,\delta)>0$ such that \eqref{eqpropM1} holds together with $T\geq T_{\mathrm{II}}$ and $\sigma \in [1/2+(\log T)^{-\theta},\sigma_1]$, where the implied constant depends only on $F$ and $\sigma_1$. 
\end{proposition}

We first prove Proposition \ref{propM1} in Section \ref{sec3.1}. 
We sometimes omit details of the proofs there since they strongly follow Guo's method in \cite{Guo1996a}. 
Towards the proof of Theorem \ref{thmM1}, we next consider in Section \ref{sec3.2} the growth of the function $\widetilde{M}_\sigma(z;F)$ of \eqref{eqpropM1}. 
We finally complete the proof of Theorem \ref{thmM1} in Section \ref{sec3.3}.

\subsection{Proof of Proposition \ref{propM1}}\label{sec3.1}
Let $F(s)$ be a function satisfying axiom $\mathrm{(4)}$. 
Then we see that 
\begin{gather*}
\frac{F'}{F}(s)
=-\sum_{n=1}^{\infty} \frac{\Lambda_F(n)}{n^s}, 
\qquad 
\sigma>1, 
\end{gather*}
where $\Lambda_F(n)$ is given by $\Lambda_F(n)=(\alpha_1(p)^m+\cdots+\alpha_g(p)^m) \log{p}$ if $n=p^m$ and $\Lambda_F(n)=0$ otherwise. 
In this section, we approximate $(F'/F)(\sigma+it)$ by some Dirichlet polynomials. 
First, we define
\begin{gather*}
w_X(n)= 
\begin{cases}
\displaystyle{1}
& \text{if $1 \leq n \leq X$}, \\
\displaystyle{ \frac{\log(X^2/n)}{\log{X}} }
& \text{if $X \leq n \leq X^2$} 
\end{cases}
\end{gather*}
for $X>1$. 
We approximate $(F'/F)(\sigma+it)$ by the following function $f_X(t,\sigma;F)$: 
\begin{gather*}
f_X(t,\sigma;F)
=-\sum_{n \leq X^2} \frac{\Lambda_F(n)}{n^{\sigma+it}} w_X(n). 
\end{gather*}

\begin{lemma}\label{lemstep1}
Let $F(s)$ be a function satisfying axioms $\mathrm{(1)}$--$\mathrm{(4)}$. 
Let $\sigma_1$ be a large fixed positive real number. 
Let $\epsilon>0$ be a small fixed real number. 
Then there exists an absolute constant $T_0>0$ such that we have 
\begin{gather}\label{eqstep1}
\frac{1}{T} \int_{0}^{T} \psi_z\left(\frac{F'}{F}(\sigma+it)\right) \,dt 
=\frac{1}{T} \int_{0}^{T} \psi_z(f_X(t,\sigma;F)) \,dt
+E_1 
\end{gather}
for all $T \geq T_0$, for all $\sigma \in [1-b^{-1}+\epsilon, \sigma_1]$, and for all $z \in \mathbb{C}$. 
The error term $E_1$ is estimated as for any $X,Y>1$ 
\begin{align}\label{eqE1}
E_1
\ll 
&\frac{1}{T} 
+Y T^{-\frac{b}{2} \{\sigma-(1-b^{-1}+\frac{\epsilon}{2})\}} \\
&+\frac{|z|}{\log{X}}
\left( \frac{X \log{Y} \log{T}}{Y} 
+\frac{X^{-\frac{1}{2} \{\sigma-(1-b^{-1}+\frac{\epsilon}{2})\}} \log{T}}{\{\sigma-(1-b^{-1}+\frac{\epsilon}{2})\}^2} 
+\frac{X}{T}
+X^{-\sigma}\log^2{T} \right), \nonumber
\end{align}
where the implied constant depends only on $F$. 
If $F(s)$ further satisfies axiom $\mathrm{(6)}$, then \eqref{eqstep1} holds with $\sigma \in [1/2+(\log T)^{-\theta}, \sigma_1]$, and we have 
\begin{align}\label{eqE1'}
E_1
\ll 
&\frac{1}{T}
+Y T^{-\frac{c}{2} (\sigma-\frac{1}{2})} (\log{T})^A \\
&+\frac{|z|}{\log{X}}
\left( \frac{X \log{Y} \log{T}}{Y} 
+\frac{X^{-\frac{1}{2} (\sigma-\frac{1}{2})} \log{T}}{(\sigma-\frac{1}{2})^2} 
+\frac{X}{T}
+X^{-\sigma}\log^2{T} \right), \nonumber
\end{align}
where the implied constant depends only on $F$. 
\end{lemma}

\begin{proof}
This lemma is an analogue of Lemma 2.1.4 of \cite{Guo1996a}. 
Let $\mathscr{B}_Y(\sigma,T;F)$ be the set of all $t \in [0,T]$ for which $|\gamma-t| \leq Y$ holds with some zeros $\rho=\beta+i \gamma$ of $F(s)$ satisfying $\beta \geq \frac{1}{2} (\sigma+1-b^{-1}+\frac{\epsilon}{2})$. 
Then we see that $E_1$ is 
\begin{gather}\label{eqE1-1}
\ll \frac{1}{T}
+\frac{\nu_1(\mathscr{B}_Y(\sigma,T;F))}{T} 
+\frac{|z|}{T} \int_{[1,T] \cap \mathscr{B}_Y(\sigma,T;F)^c}
\left|\frac{F'}{F}(\sigma+it)-f_X(t,\sigma;F)\right| \,dt, 
\end{gather}
since $|\psi_z(w)-\psi_z(w')| \leq |z||w-w'|$. 
By the definition of $\mathscr{B}_Y(\sigma,T;F)$, we have 
\begin{gather*}
\nu_1(\mathscr{B}_Y(\sigma,T;F))
\leq 2Y N_F\left(\frac{1}{2} \left(\sigma+1-b^{-1}+\frac{\epsilon}{2}\right), T\right). 
\end{gather*}
Furthermore, estimate \eqref{eq:WZD} implies that the second term of \eqref{eqE1-1} is 
\begin{gather*}
\ll Y T^{1-\frac{b}{2} \{\sigma-(1-b^{-1}+\frac{\epsilon}{2})\}} 
\end{gather*}
for $\sigma>1-b^{-1}+\epsilon/2$. 
Then we estimate the third term. 
For this, Guo used the formula of \cite[Lemma 2]{Selberg1943}, and we need a similar formula for general $F(s)$. 
We first recall that the following estimate 
\begin{gather}\label{eqEsti}
\frac{F'}{F}(s) 
\ll \log^2{(|t|+2)} 
\end{gather}
holds if $s=\sigma+it$ satisfies $-1 \leq \sigma \leq 2$ and has distance $\gg \log(|t|+2)^{-1}$ from zeros and poles of $F(s)$. 
This can be easily deduced from axioms  $\mathrm{(1)}$--$\mathrm{(4)}$. 
Let $c=\max\{2,1+\sigma\}$ and choose $T_m \in (m, m+1]$ and $0<\delta<1$ such that the edges $[c+i T_m, -\delta+i T_m]$, $[c-i T_m, -\delta-i T_m]$, and $[-\delta-i T_m, -\delta+i T_m]$ have distance $\gg \log(|t|+2)^{-1}$ from zeros and poles of $F(s)$. 
Then, we consider the integral 
\begin{gather*}
\frac{1}{2\pi i} \int_{c-i T_m}^{c+i T_m} \frac{F'}{F}(z) \frac{X^{z-s}-X^{2(z-s)}}{(z-s)^2} \,dz. 
\end{gather*}
We see that 
\begin{gather*}
\lim_{m \to\infty} \frac{1}{2\pi i} \int_{c-i T_m}^{c+i T_m} \frac{F'}{F}(z) \frac{X^{z-s}-X^{2(z-s)}}{(z-s)^2} \,dz 
=-f_X(t,\sigma;F) \log{X} 
\end{gather*}
and change the contour by the edges $[c+i T_m, -\delta+i T_m]$, $[c-i T_m, -\delta-i T_m]$, and $[-\delta-i T_m, -\delta+i T_m]$. 
The integrals on the horizontal edges tend to $0$ as $m \to\infty$ due to estimate \eqref{eqEsti}, and we have also by \eqref{eqEsti}, 
\begin{gather*}
\frac{1}{2\pi i} \int_{-\delta-i T_m}^{-\delta+i T_m} \frac{F'}{F}(z) \frac{X^{z-s}-X^{2(z-s)}}{(z-s)^2} \,dz 
\ll_{\sigma_0} X^{-\sigma} \log^2{T} 
\end{gather*}
for any $\sigma \geq \sigma_0>0$ and $t \in [1,T]$. 
Calculating the residues, we obtain the following formula: 
\begin{align}\label{eqSelb}
\frac{F'}{F}(s)
&=f_X(t,\sigma;F) 
-\frac{m_1}{\log{X}} \frac{X^{1-s}-X^{2(1-s)}}{(1-s)^2} 
+\frac{m_0}{\log{X}} \frac{X^{-s}-X^{-2s}}{s^2} \\
&+\frac{1}{\log{X}} \sum_{\rho} \frac{X^{\rho-s}-X^{2(\rho-s)}}{(\rho-s)^2} 
+O_{\sigma_0}\left( \frac{1}{\log{X}} X^{-\sigma} \log^2{T} \right), \nonumber
\end{align}
where $m_1, m_0 \geq0$ are orders of the possible pole of $F(s)$ at $s=1$ and the possible zero of $F(s)$ at $s=0$, respectively, and $\rho$ runs through nontrivial zeros of $F(s)$. 
In order to complete the proof of Lemma \ref{lemstep1}, we must consider the contributions of the second, third, and fourth terms of \eqref{eqSelb}. 
They are estimated by an argument similar to the proof of Lemma 2.1.4 of \cite{Guo1996a}. 
Thus we find the first part of Lemma \ref{lemstep1}. 

All changes that we need for the proof of the second part are just replacing the definition of $\mathscr{B}_Y(\sigma,T;F)$ with the set of all $t \in [0,T]$ for which $|\gamma-t| \leq{Y}$ holds with some zeros $\rho=\beta+i \gamma$ of $F(s)$ satisfying $\beta \geq \frac{1}{2}(\sigma+\frac{1}{2})$. 
By the axiom $\mathrm{(6)}$, we have 
\begin{gather*}
\nu_1(\mathscr{B}_Y(\sigma,T;F))
\leq 2Y N_F\left(\frac{1}{2}\left(\sigma+\frac{1}{2}\right), T \right) 
\ll YT^{1-\frac{c}{2} (\sigma-\frac{1}{2})} (\log{T})^A. 
\end{gather*}
The remaining estimates are given in a similar way. 
\end{proof}

Towards the next step, we define 
\begin{gather*}
g_X(t,\sigma;F)
=-\sum_{n \leq X^2} \frac{\Lambda_F(n)}{n^{\sigma+it}} 
\quad\text{and}\quad
h_X(t,\sigma;F)
=-\sum_{p \leq X^2} \sum_{m=1}^{\infty} \frac{\Lambda_F(p^m)}{p^{m(\sigma+it)}}
\end{gather*}
for $X>1$. 
Then we have the following three lemmas: 

\begin{lemma}\label{lemstep2}
Let $F(s)$ be a function satisfying axioms $\mathrm{(1)}$ and $\mathrm{(4)}$. 
Then there exists an absolute constant $T_0>0$ such that we have 
\begin{gather*}
\frac{1}{T} \int_{0}^{T} \psi_z(f_X(t,\sigma;F)) \,dt 
=\frac{1}{T} \int_{0}^{T} \psi_z(g_X(t,\sigma;F)) \,dt
+E_2
\end{gather*}
for all $T \geq T_0$, for all $\sigma>1/2$, and for all $z \in \mathbb{C}$. 
The error term $E_2$ is estimated as 
\begin{gather}\label{eqE2}
E_2
\ll \frac{g |z| \log{X}}{(2\sigma-1)^{\frac{1}{2}}} \left(1+\frac{X^2}{T}\right)^{\frac{1}{2}} X^{\frac{1}{2}-\sigma} 
\end{gather}
for any $X>1$. 
The implied constant is absolute. 
\end{lemma}

\begin{lemma}\label{lemstep3}
Let $F(s)$ be a function satisfying axioms $\mathrm{(1)}$ and $\mathrm{(4)}$. 
Then there exists an absolute constant $T_0>0$ such that we have 
\begin{gather*}
\frac{1}{T} \int_{0}^{T} \psi_z(g_X(t,\sigma;F)) \,dt 
=\frac{1}{R} \int_{0}^{R} \psi_z(g_X(r,\sigma;F)) \,dr
+E_3
\end{gather*}
for all $R \geq T \geq T_0$, for all $\sigma>1/2$, and for all $z \in \mathbb{C}$. 
The error term $E_3$ is estimated as 
\begin{align}\label{eqE3}
E_3
&\ll \frac{ g^N X^{5N} }{T} (1+|z|^2)^{\frac{N}{2}} \\
&\qquad+\frac{(8g|z|)^N}{N!} \left(1+\frac{X^N}{T}\right) 
\left\{ (\zeta(2\sigma)^{\frac{1}{2}} \log{X})^N \left(\frac{N}{2}\right)! 
+\zeta'(2\sigma)^N \right\} \nonumber
\end{align}
for any $X>1$ and any large even integer $N$. 
The implied constant is absolute. 
\end{lemma}

\begin{lemma}\label{lemstep4}
Let $F(s)$ be a function satisfying axioms $\mathrm{(1)}$ and $\mathrm{(4)}$. 
Then there exists an absolute constant $T_0>0$ such that we have 
\begin{gather*}
\frac{1}{R} \int_{0}^{R} \psi_z(g_X(r,\sigma;F)) \,dr 
=\frac{1}{R} \int_{0}^{R} \psi_z(h_X(r,\sigma;F)) \,dr
+E_4
\end{gather*}
for all $R \geq T \geq T_0$, for all $\sigma>1/2$, and for all $z \in \mathbb{C}$. 
The error term $E_4$ is estimated as 
\begin{gather}\label{eqE4}
E_4
\ll \frac{g |z| \log{X}}{2\sigma-1} X^{1-2\sigma} 
\end{gather}
for any $X>1$. 
The implied constant is absolute. 
\end{lemma}

These lemmas are analogues of Lemmas 2.2.5, 2.1.6, and 2.1.10 in \cite{Guo1996a}. 
Note that we have $|\Lambda_F(n)| \leq g \Lambda(n)$ due to axioms $\mathrm{(1)}$ and $\mathrm{(4)}$, where $\Lambda(n)=\Lambda_\zeta(s)$ is the usual von Mangolt function. 
In fact, by axiom $\mathrm{(4)}$ we have $\Lambda_F(p^m)=(\alpha_1(p)^m+\ldots+\alpha_g(p)^m) \log{p}$, and by axiom $\mathrm{(1)}$ the absolute values of $\alpha_j(p)$ are less than or equal to $1$; see Lemma 2.2 of \cite{Steuding2007}. 
Therefore we obtain these lemmas by replacing $\Lambda(n)$ with $\Lambda_F(n)$ in the proofs of the corresponding lemmas in \cite{Guo1996a}. 

Let $F(s)$ be a function satisfying axioms $\mathrm{(1)}$--$\mathrm{(4)}$. 
Let $\sigma_1$ be a large fixed positive real number. 
Let $\epsilon>0$ be a small fixed real number. 
By the above lemmas, we have for all $R \geq T \geq T_0$ and for all $\sigma \in [1-b^{-1}+\epsilon, \sigma_1]$, 
\begin{gather}\label{eqstep1234}
\frac{1}{T} \int_{0}^{T} \psi_z\left(\frac{F'}{F}(\sigma+it)\right) \,dt 
=\frac{1}{R} \int_{0}^{R} \psi_z(h_X(r,\sigma;F)) \,dr
+E_1+E_2+E_3+E_4, 
\end{gather}
where the error terms $E_j$ are estimated as in \eqref{eqE1}, \eqref{eqE2}, \eqref{eqE3}, and \eqref{eqE4}. 
Let $\theta, \delta>0$ with $\delta+3\theta<1/2$. 
We take $X$, $Y$, and $N$ as the following functions in $T$: 
\begin{gather*}
X
=\exp((\log{T})^{\theta_1}), 
\quad
Y
=\exp((\log{T})^{\theta_2}), 
\quad\text{and}\quad
N
=2\lfloor(\log{T})^{\theta_3}\rfloor,
\end{gather*}
where $\theta_1=(5/3)\theta$, $\theta_2=(\theta_1+1-\theta)/2$, $\theta_3=((2\delta+\theta+2\theta_1)+(1-\theta_1))/2$. 
Moreover, let $T'_0=T'_0(\theta, \epsilon) \geq T_0$ with 
\begin{gather*}
(\log T'_0)^{-\theta}
\leq \epsilon/2. 
\end{gather*}
Then we have $\sigma \geq 1-b^{-1}+\epsilon/2+(\log T)^{-\theta}$ for $T \geq T'_0$. 
Hence, there exists a positive real number $T_{\mathrm{I}}=T_{\mathrm{I}}(F,\theta,\delta,\epsilon) \geq T'_0$ such that we have 
\begin{gather}\label{eqE1234}
E_1+E_2+E_3+E_4
\ll\exp\left(-\frac{1}{4} (\log T)^{\frac{2}{3}\theta}\right)
\end{gather}
for all $T \geq T_{\mathrm{I}}$ and for all $z \in \Omega$ with the implied constant depending only on $F$ and $\epsilon$. 

Then, let $F(s)$ further satisfy axiom $\mathrm{(6)}$. 
In this case, we obtain that the formula \eqref{eqstep1234} holds for all $R \geq T \geq T_0$ and for all $\sigma \in [1/2+(\log{T})^{-\theta}, \sigma_1]$, where the error terms $E_j$ are estimated as in \eqref{eqE1'}, \eqref{eqE2}, \eqref{eqE3}, and \eqref{eqE4}. 
Therefore there exists a positive real number $T_{\mathrm{II}}=T_{\mathrm{II}}(F,\theta,\delta)>T_0$ such that we have the same estimate as \eqref{eqE1234} for all $T \geq T_{\mathrm{II}}$ and for all $z \in \Omega$. 

Next, applying Lemma 2 of \cite{HeathBrown1992}, we see that 
\begin{gather}\label{eqstep5}
\lim_{R \to\infty} \frac{1}{R} \int_{0}^{R}\psi_z(h_X(r,\sigma;F)) \,dr 
=\prod_{p \leq X^2} \int_{0}^{1} 
\psi_z\left( \sum_{m=1}^{\infty} \frac{\Lambda_F(p^m)}{p^{m\sigma}} e^{2 \pi i m \theta} \right) \,d\theta 
\end{gather}
since the system 
\begin{gather*}
\left\{ \frac{\log p}{2\pi} ~\middle|~ \text{$p$ is a prime number} \right\}
\end{gather*}
is linearly independent over $\mathbb{Q}$. 
We define 
\begin{gather}\label{eqMp}
\widetilde{M}_{\sigma,p}(z;F)
=\int_{0}^{1} \psi_z\left( \sum_{m=1}^{\infty} \frac{\Lambda_F(p^m)}{p^{m\sigma}} e^{2 \pi i m \theta} \right) \,d\theta. 
\end{gather}
Then we obtain the following lemma on $\widetilde{M}_{\sigma,p}(z;F)$, which is proved in Section \ref{sec3.2}. 

\begin{lemma}\label{lemstep5}
Let $F(s)$ be a function satisfying axioms $\mathrm{(1)}$ and $\mathrm{(4)}$. 
Let $\sigma_1$ be a large fixed positive real number. 
Let $\theta, \delta>0$ be real numbers with $\delta+3\theta<1/2$. 
Then there exists a positive real number $T_0=T_0(F,\sigma_1,\theta,\delta)$ such that we have 
\begin{gather*}
\prod_{p>X^2} \widetilde{M}_{\sigma,p}(z;F) 
=1+O\left( \exp\left(-\frac{1}{4} (\log T)^{\frac{2}{3}\theta}\right) \right)
\end{gather*}
for all $T \geq T_0$, for all $\sigma \in [1/2+(\log T)^{-\theta}, \sigma_1]$, and for all $z \in \Omega$. 
Here we denote $X=\exp((\log T)^{\frac{5}{3}\theta})$, and the implied constant depends only on $F$ and $\sigma_1$. 
\end{lemma}

We prove Proposition \ref{propM1} with the above preliminary lemmas. 

\begin{proof}[Proof of Proposition \ref{propM1}]
By \eqref{eqstep1234}, \eqref{eqE1234}, and \eqref{eqstep5}, we have 
\begin{gather*}
\frac{1}{T} \int_{0}^{T} \psi_z\left(\frac{F'}{F}(\sigma+it)\right) \,dt 
=\prod_{p \leq X^2} \widetilde{M}_{\sigma,p}(z;F) 
+O\left( \exp\left(-\frac{1}{4} (\log T)^{\frac{2}{3}\theta}\right) \right). 
\end{gather*}
We consider the replacement of the product $\prod_{p \leq X^2} \widetilde{M}_{\sigma,p}(z;F)$ with $\prod_{p} \widetilde{M}_{\sigma,p}(z;F)$, where the error is estimated as 
\begin{gather*}
\left| \prod_{p} \widetilde{M}_{\sigma,p}(z;F)-\prod_{p \leq X^2} \widetilde{M}_{\sigma,p}(z;F) \right| 
\leq \left| \prod_{p>X^2} \widetilde{M}_{\sigma,p}(z;F)-1 \right|, 
\end{gather*}
since $|\widetilde{M}_{\sigma,p}(z;F)| \leq 1$ by definition. 
Hence we have 
\begin{gather*}
\prod_{p \leq X^2} \widetilde{M}_{\sigma,p}(z;F) 
=\prod_{p} \widetilde{M}_{\sigma,p}(z;F) 
+O\left(\exp\left(-\frac{1}{4} (\log T)^{\frac{2}{3}\theta}\right) \right) 
\end{gather*}
by Lemma \ref{lemstep5}. 
Therefore Proposition \ref{propM1} follows if we define 
\begin{gather*}
\widetilde{M}_\sigma(z;F)
=\prod_{p} \widetilde{M}_{\sigma,p}(z;F). 
\end{gather*}
\end{proof}

\subsection{Estimates on $\widetilde{M}_\sigma(z;F)$}\label{sec3.2}
In this section, we examine some analytic properties of the function $\widetilde{M}_\sigma(z;F)$. 
By definition \eqref{eqMp} and $\psi_z(w)=\exp(i \RE(z \overline{w}))$, we have 
\begin{gather*}
\widetilde{M}_{\sigma,p}(z;F) 
=\int_{0}^{1} \exp(i x a_p(\theta,\sigma;F) + i y b_p(\theta,\sigma;F)) \,d\theta, 
\end{gather*}
where $z=x+i y$ and $a_p(\theta,\sigma;F), b_p(\theta,\sigma;F)$ are functions such that 
\begin{align*}
a_p(\theta,\sigma;F) 
&=\sum_{m=1}^{\infty} \frac{1}{p^{m\sigma}} 
\{\RE \Lambda_F(p^m) \cos(2 \pi m \theta)-\IM \Lambda_F(p^m) \sin(2 \pi m \theta)\}, \\
b_p(\theta,\sigma;F) 
&=\sum_{m=1}^{\infty} \frac{1}{p^{m\sigma}}
\{\RE \Lambda_F(p^m) \sin(2 \pi m \theta)+\IM \Lambda_F(p^m) \cos(2 \pi m \theta)\}. 
\end{align*}
Then we define 
\begin{gather}\label{eqMp'}
\widetilde{M}_p(s,z_1,z_2;F) 
=\int_{0}^{1} \exp(i z_1 a_p(\theta,s;F)+i z_2 b_p(\theta,s;F)) \,d\theta
\end{gather}
for $\RE s>0$ and $z_1, z_2 \in \mathbb{C}$. 
We have $\widetilde{M}_{\sigma,p}(x+iy;F)=\widetilde{M}_p(\sigma,x,y;F)$ if $\sigma>0$ and $x,y \in \mathbb{R}$. 
For the study on the function $\widetilde{M}_p(s,z_1,z_2;F)$, the following lemma is fundamental, which is easily deduced from the expansion of $\exp(z)$ and the calculations of integrals. 

\begin{lemma}\label{lemFund}
Let $F(s)$ be a function that satisfies axiom $\mathrm{(4)}$. 
Then we have 
\begin{gather}\label{eqFund} 
\widetilde{M}_p(s,z_1,z_2;F) 
=1- \mu_p + R_p 
\end{gather}
for $\sigma=\RE s>0$ and $z_1, z_2 \in \mathbb{C}$, where 
\begin{align*}
\mu_p 
&=\mu_p(s,z_1,z_2;F) 
=\frac{z_1^2+z_2^2}{4} \sum_{m=1}^{\infty} \frac{|\Lambda_F(p^m)|^2}{p^{2ms}}, \\
R_p 
&=R_p(s,z_1,z_2;F) 
=\int_{0}^{1} \sum_{k=3}^{\infty} \frac{i^k}{k!} \{z_1 a_p(\theta,s;F) + z_2 b_p(\theta,s;F)\}^k \,d\theta.
\end{align*}
\end{lemma}

Therefore, if $\mu_p$ and $R_p$ are sufficiently small, we have 
\begin{gather}\label{eqlogMp}
\log\widetilde{M}_p(s,z_1,z_2;F)
=-\mu_p+R_p+O(|\mu_p|^2+|R_p|^2), 
\end{gather}
where $\log$ is the principal blanch of logarithm. 
Using Lemma \ref{lemFund}, we study the function 
\begin{gather}\label{eqprod}
\widetilde{M}(s,z_1,z_2;F)
=\prod_{p} \widetilde{M}_p(s,z_1,z_2;F). 
\end{gather}

\begin{proposition}\label{propM2}
Let $F(s)$ be a function satisfying axioms $\mathrm{(1)}$ and $\mathrm{(4)}$. 
Assume that $(s, z_1, z_2)$ varies on $\{\RE s>1/2\} \times \mathbb{C} \times \mathbb{C}$. 
If we fix two of the variables, the function $\widetilde{M}(s,z_1,z_2;F)$ is holomorphic with respect to the reminder variable. 
\end{proposition}

\begin{proof}
Let $K$ be any compact subset on the half plane $\{\RE s>1/2\}$, and let $K_1, K_2$ be any compact subsets on $\mathbb{C}$. 
Assume that $(s, z_1, z_2) \in K \times K_1 \times K_2$, and let $\sigma_0$ be the smallest  real part of $s \in K$. 
As in Section \ref{sec3.1}, we have $|\Lambda(p^m)| \leq g \log{p}$, where $g$ is the constant in axiom $\mathrm{(4)}$. 
Then we obtain 
\begin{gather*}
\mu_p
\ll \frac{g^2 (\log{p})^2}{p^{2\sigma_0}}
\quad\text{and}\quad
R_p
\ll \frac{g^3 (\log{p})^3}{p^{3\sigma_0}}, 
\end{gather*}
where the implied constants depend only on $K, K_1, K_2$. 
Thus, by \eqref{eqlogMp}, we have $\log\widetilde{M}_p(s,z_1,z_2;F) \ll g^2 (\log{p})^2 p^{-2\sigma_0}$ for all $p>M$, where $M=M(K,K_1,K_2)$ is a sufficiently large constant that depends only on $K$, $K_1$, and $K_2$. 
The series $\sum_{p} (\log{p})^2 p^{-2\sigma_0}$ converges since $\sigma_0>1/2$; therefore infinite product \eqref{eqprod} uniformly converges on $K \times K_1 \times K_2$. 
Every local parts $\widetilde{M}_p(s,z_1,z_2;F)$ are holomorphic, and hence we have the result. 
\end{proof}

We estimate the growth of $\widetilde{M}(s,z_1,z_2;F)$ with $z_1$ and $z_2$ near the real axis. 

\begin{proposition}\label{propM3}
Let $F(s)$ be a function satisfying axioms $\mathrm{(1)}$, $\mathrm{(4)}$, and $\mathrm{(5)}$. 
Let $\sigma>1/2$ be an arbitrarily fixed real number. 
Then there exist positive constants $K=K(\sigma;F)$ and $c=c(\sigma;F)$ such that for all $x, y \in \mathbb{R}$ with $|x|+|y| \geq K$, and for all non-negative integers $m$ and $n$, we have 
\begin{gather*}
\frac{\partial^{m+n}}{\partial z_1^m \partial z_2^n} 
\widetilde{M}(\sigma,z_1,z_2;F) 
\ll \exp\left( -c (|x|+|y|)^{\frac{1}{\sigma}} (\log(|x|+|y|))^{\frac{1}{\sigma}-1} \right) 
\end{gather*}
for any $z_1, z_2 \in \mathbb{C}$ with $|z_1-x|<1/4$, $|z_2-y|<1/4$. 
The implied constant depends only on $m$ and $n$. 
\end{proposition}

\begin{proof}
Let $K>1$ and $c_0<1$ be positive constants chosen later, and assume that $x,y \in \mathbb{R}$ with $|x|+|y| \geq K$. 
We define 
\begin{gather*}
P_0 
=\left( \frac{g(|x|+|y|)}{c_0} \log\frac{g(|x|+|y|)}{c_0} \right)^{\frac{1}{\sigma}} 
\end{gather*}
for any fixed $\sigma>1/2$. 
Then for any $p \geq P_0$, we see that 
\begin{gather*}
\frac{(|x|+|y|) g \log{p}}{p^\sigma} 
\leq \frac{(|x|+|y|) g \log{P_0}}{P_0^\sigma} 
\leq c_0 c_1
\end{gather*}
with an absolute constant $c_1>0$. 
Hence, we estimate $\mu_p$ and $R_p$ in Lemma \ref{lemFund} arbitrarily small if we let the constant $c_0$ suitably small. 
Thus formula \eqref{eqlogMp} holds. 
We then replace $\mu_p$ in \eqref{eqlogMp} with the real number 
\begin{gather*}
\mu'_p 
=\mu'_p(\sigma,x,y;F) 
=\frac{x^2+y^2}{4} \sum_{m=1}^{\infty} \frac{|\Lambda_F(p^m)|^2}{p^{2m \sigma}}. 
\end{gather*}
The error of the replacement is estimated as 
\begin{gather*}
|\mu_p-\mu'_p|
\leq (|x|+|y|) \sum_{m=1}^{\infty} \frac{|\Lambda_F(p^m)|^2}{p^{2m \sigma}}
\end{gather*}
if we assume that $|z_1-x|<1/2$ and $|z_2-y|<1/2$. 
Moreover, we have 
\begin{align*}
\mu_p^2
&\ll \left(\frac{(|x|+|y|) g \log{p}}{p^\sigma}\right)^4 
\leq \frac{(|x|+|y|) g \log{P_0}}{P_0^\sigma} \left(\frac{(|x|+|y|) g \log{p}}{p^\sigma}\right)^3 \\
&\ll \frac{(|x|+|y|)^3 g^3 (\log{p})^3}{p^{3\sigma}}
\end{align*}
and 
\begin{gather*}
R_p 
\ll (|x|+|y|)^3 \left(\sum_{m=1}^{\infty} \frac{|\Lambda_F(p^m)|}{p^{m \sigma}}\right)^3 
\ll \frac{(|x|+|y|)^3 g^3 (\log{p})^3}{p^{3\sigma}}, 
\end{gather*}
where all implied constants are absolute. 
Therefore by \eqref{eqlogMp} we have for any $p \geq P_0$, 
\begin{align*}
&\left| \log\widetilde{M}_p(\sigma,z_1,z_2;F) 
+\frac{x^2+y^2}{4} \sum_{m=1}^{\infty} \frac{|\Lambda_F(p^m)|^2}{p^{2m \sigma}} \right| \\
&\leq (|x|+|y|) \sum_{m=1}^{\infty} \frac{|\Lambda_F(p^m)|^2}{p^{2m \sigma}} 
+B (|x|+|y|)^3 \frac{g^3 (\log{p})^3}{p^{3\sigma}} 
\end{align*}
with some absolute constant $B>0$. 
Thus for sufficiently large $K$, if $|x|+|y| \geq K$, then we obtain 
\begin{align*}
&\RE \log\widetilde{M}_p(\sigma,z_1,z_2;F) \\
&\leq -A (|x|+|y|)^2 \sum_{m=1}^{\infty} \frac{|\Lambda_F(p^m)|^2}{p^{2m \sigma}} 
+B (|x|+|y|)^3 \frac{g^3 (\log{p})^3}{p^{3\sigma}} \\
&\leq -A (|x|+|y|)^2 \frac{|\Lambda_F(p)|^2}{p^{2\sigma}} 
+B (|x|+|y|)^3 \frac{g^3 (\log{p})^3}{p^{3\sigma}} 
\end{align*}
with some absolute constant $A>0$. 
Note that 
\begin{gather*}
\Lambda_F(p) 
=(\alpha_1(p)+\cdots\alpha_g(p)) \log{p} 
=-a_F(p) \log{p}
\end{gather*}
from axiom $\mathrm{(4)}$. 
Hence, we have 
\begin{align}\label{eqlarge1}
&\left| \prod_{p \geq P_0} \widetilde{M}_p(\sigma,z_1,z_2;F) \right| \\
&\leq \exp\left(-A (|x|+|y|)^2 \sum_{p \geq P_0} \frac{|a_F(p)|^2 (\log{p})^2}{p^{2\sigma}} 
+B (|x|+|y|)^3 g^3 \sum_{p \geq P_0} \frac{(\log{p})^3}{p^{3\sigma}}\right) \nonumber\\
&\leq \exp\left(-A (|x|+|y|)^2 \sum_{p \geq P_0} \frac{(\log{p})^2}{p^{2\sigma}} |a_F(p)|^2 
+B c_0 c_1 g^2 (|x|+|y|)^2 \sum_{p \geq P_0} \frac{(\log{p})^2}{p^{2\sigma}}\right). \nonumber
\end{align}
Then we estimate 
\begin{gather*}
\sum_{p \geq P_0} \frac{(\log{p})^2}{p^{2\sigma}} |a_F(p)|^2 
\quad\text{and}\quad
\sum_{p \geq P_0} \frac{(\log{p})^2}{p^{2\sigma}}. 
\end{gather*}
We see that for any $\sigma>1/2$, there exists a constant $X_0(\sigma; F)>0$ such that for any $X \geq X_0(\sigma; F)$, 
\begin{align*}
\sum_{p \geq X} \frac{(\log{p})^2}{p^{2\sigma}} |a_F(p)|^2 
&\geq \frac{\kappa}{2(2\sigma-1)} X^{1-2\sigma} \log{X} \\
\sum_{p \geq X} \frac{(\log{p})^2}{p^{2\sigma}}
&\leq \frac{2}{2\sigma-1} X^{1-2\sigma} \log{X}.
\end{align*}
Indeed, the first inequality is deduced by summing by parts with axiom $\mathrm{(5)}$, and we obtain the second inequality in a similar way. 
Then, we let $c_0=c_0(F)$ smaller so that $2B c_0 c_1 g^2<A \kappa/2$. 
If we let $K=K(\sigma; F)$ suitably large, then we obtain for $|u|+|v| \geq K$, 
\begin{align*}
&-A (|x|+|y|)^2 \sum_{p \geq P_0} \frac{(\log{p})^2}{p^{2\sigma}} |a_F(p)|^2 
+B c_0 c_1 g^2 (|x|+|y|)^2 \sum_{p \geq P_0} \frac{(\log{p})^2}{p^{2\sigma}} \\
&\leq -c (|x|+|y|)^{\frac{1}{\sigma}} (\log(|x|+|y|))^{\frac{1}{\sigma}-1} 
\end{align*}
with some positive constant $c=c(\sigma;F)$. 
Hence we obtain 
\begin{gather}\label{eqlarge2}
\left| \prod_{p \geq P_0} \widetilde{M}_p(\sigma,z_1,z_2;F) \right| 
\leq \exp\left(-c (|x|+|y|)^{\frac{1}{\sigma}} (\log(|x|+|y|))^{\frac{1}{\sigma}-1} \right). 
\end{gather}
The estimate on the contributions of $\widetilde{M}_p(\sigma,z_1,z_2;F)$ for $p<P_0$ remains. 
By definition \eqref{eqMp'}, we see that 
\begin{align*}
\left|\widetilde{M}_p(\sigma,z_1,z_2;F)\right| 
&\leq \int_{0}^{1} \exp(-\IM(z_1) a_p(\theta,\sigma;F)-\IM(z_2) b_p(\theta,\sigma;F)) \,d\theta \\
&\leq \int_{0}^{1} \exp( |a_p(\theta,\sigma;F)|+|b_p(\theta,\sigma;F)| ) \,d\theta \\
&\leq \exp\left( C \frac{g \log{p}}{p^\sigma} \right) 
\end{align*}
with some absolute positive constant $C$ since $|z_1-x|<1/2$, $|z_2-y|<1/2$, and $x,y \in \mathbb{R}$. 
Thus we have 
\begin{gather*}
\left| \prod_{p<P_0} \widetilde{M}_p(\sigma,z_1,z_2;F) \right| 
\leq \exp\left( C \sum_{p<P_0} \frac{g\log{p}}{p^\sigma} \right) 
\leq \exp\left( C g \log{P_0} P_0^{\frac{1}{2}} \right). 
\end{gather*}
Then we see that for $|x|+|y| \geq K$, 
\begin{gather}\label{eqsmall}
\left| \prod_{p<P_0} \widetilde{M}_p(\sigma,z_1,z_2;F) \right| 
\leq\exp\left( C' (|x|+|y|)^{\frac{3}{4\sigma}} \right), 
\end{gather}
where $C'=C'(F)$ is some positive constant. 
Therefore we obtain
\begin{gather}\label{eqMEsti}
\left|\widetilde{M}(\sigma,z_1,z_2;F)\right| 
\leq\exp\left(-c (|x|+|y|)^{\frac{1}{\sigma}} (\log(|x|+|y|))^{\frac{1}{\sigma}-1} \right) 
\end{gather}
by \eqref{eqlarge2} and \eqref{eqsmall}, where $c=c(\sigma; F)$ is some positive constant. 
We finally assume that $|z_1-x|<1/4$ and $|z_2-y|<1/4$. 
Then, applying Cauchy's integral formula, we have 
\begin{gather*}
\frac{\partial^{m+n}}{\partial z_1^m \partial z_2^n} 
\widetilde{M}(\sigma,z_1,z_2;F) 
=\frac{m!n!}{(2\pi i)^2} \iint_{\substack{|\xi_1-z_1|=1/4, \\ |\xi_2-z_2|=1/4}}
\frac{\widetilde{M}(\sigma,\xi_1,\xi_2;F)}{(\xi_1-z_1)^{m+1} (\xi_2-z_2)^{n+1}} \,d \xi_1 d \xi_2. 
\end{gather*}
Therefore by estimate \eqref{eqMEsti}, the desired result follows. 
\end{proof}

\begin{remark}\label{remM}
We find that $\widetilde{M}_\sigma(z;F)$ is a Schwartz function according to Proposition \ref{propM3}. 
Hence its Fourier inverse 
\begin{gather*}
M_\sigma(z;F) 
=\int_{\mathbb{C}} \widetilde{M}_\sigma(w;F) \psi_{-z}(w) \,|dw| 
\end{gather*}
is also a Schwartz function, and belongs to the class $\Lambda$. 
Thus we have $\widetilde{M}_\sigma(z;F)=(M_\sigma(z;F))^{\wedge}$. 
By a simple calculation, we see that $M_\sigma(z;F)$ is real valued. 
\end{remark}

Finally, we prove Lemma \ref{lemstep5} in Section \ref{sec3.1}. 

\begin{proof}[Proof of Lemma \ref{lemstep5}]
Assume $p \geq X^2$ with $X=\exp((\log T)^{\frac{5}{3}\theta})$. 
Then we see that $\mu_p=\mu_p(\sigma,x,y;F)$ and $R_p=R_p(\sigma,x,y;F)$ in Lemma \ref{lemFund} are small when $T$ is sufficiently large. 
In fact, we have for $p \geq X^2$, 
\begin{gather*}
\mu_p
\ll (x^2+y^2) \frac{g^2 (\log{p})^2}{p^{2\sigma}} 
\ll \{(|x|+|y|) g X^{1-2\sigma} \log{X}\}^2. 
\end{gather*}
By the setting for $X, z=x+iy$, and $\sigma$, we have 
\begin{gather*}
X^{1-2\sigma} \log{X}
\leq \exp\left(-\frac{1}{4} (\log{T})^{\frac{2}{3}\theta} \right)
\to 0
\end{gather*}
as $T \to\infty$. 
The argument for $R$ is similar. 
Hence by \eqref{eqlogMp}, we obtain 
\begin{gather*}
\log\widetilde{M}_{\sigma,p}(z;F)
\ll (x^2+y^2) \frac{(\log{p})^2}{p^{2\sigma}}, 
\end{gather*}
where the implied constant depends only on $F$. 
Therefore we have 
\begin{align*}
\prod_{p \geq X^2} 
\widetilde{M}_{\sigma,p}(z;F) 
&=\exp\left( \sum_{p \geq X^2} \log\widetilde{M}_{\sigma,p}(z;F) \right) \\
&=1+O\left( (x^2+y^2) \sum_{p \geq X^2} \frac{(\log{p})^2}{p^{2\sigma}} \right). 
\end{align*}
Applying the prime number theorem, we estimate the above error term as 
\begin{gather*}
(x^2+y^2) \sum_{p \geq X^2} \frac{(\log{p})^2}{p^{2\sigma}} 
\ll (x^2+y^2) \frac{X^{2(1-2\sigma)} \log{X}}{(\sigma-\frac{1}{2})^2} 
\leq \exp\left(-\frac{1}{4} (\log T)^{\frac{2}{3}\theta}\right)
\end{gather*}
by the assumptions on $X$, $z=x+iy$, and $\sigma$. 
Here the implied constant depends only on $F$ and $\sigma_1$. 
\end{proof}

\subsection{Completion of the proof}\label{sec3.3}

\begin{proof}[Proof of Theorem \ref{thmM1}]
We only consider the case of $F \in \mathcal{S}_{\mathrm{I}}$ since the case $F \in \mathcal{S}_{\mathrm{II}}$ follows completely in an analogous way. 
By the definition of the class $\Lambda$, for any $\Phi \in \Lambda$ we have 
\begin{gather*}
\Phi(w) 
=\int_{\mathbb{C}} \widehat{\Phi}(z) \psi_{-z}(w) \,|dw|. 
\end{gather*}
Hence, by Proposition \ref{propM1}, we see that for all $T \geq T_I$, 
\begin{align*}
\frac{1}{T} \int_{0}^{T} \Phi\left(\frac{F'}{F}(\sigma+it)\right) \,dt 
&=\int_{\Omega} \widehat{\Phi}(z) \frac{1}{T} \int_{0}^{T} \psi_{-z}\left(\frac{F'}{F}(\sigma+it)\right) \,dt \,|dz| 
+E_1 \\
&=\int_{\Omega} \widehat{\Phi}(z) \widetilde{M}_\sigma(-z;F) \,|dz|
+E_1+E_2 \\
&=\int_{\mathbb{C}} \widehat{\Phi}(z) \widetilde{M}_\sigma(-z;F) \,|dz|
+E_1+E_2+E_3, 
\end{align*}
where the error terms are estimated as 
\begin{align*}
E_1 
&=\int_{\mathbb{C} \setminus \Omega} \widehat{\Phi}(z) 
\frac{1}{T} \int_{0}^{T} \psi_{-z}\left(\frac{F'}{F}(\sigma+it)\right) \,dt \,|dz| 
\ll \int_{\mathbb{C} \setminus \Omega} |\widehat{\Phi}(z)| \,|dz|, \\
E_2 
&\ll \exp\left(-\frac{1}{4} (\log T)^{\frac{2}{3}\theta}\right) \int_{\Omega} |\widehat{\Phi}(z)| \,|dz|, \\
E_3 
&\ll \int_{\mathbb{C} \setminus \Omega} |\widehat{\Phi}(z)| \,|dz|. 
\end{align*}
Here all implied constants depend at most only on $F$, $\sigma_1$, and $\epsilon$. 
We find that 
\begin{gather*}
\int_{\mathbb{C}} \widehat{\Phi}(w) \widetilde{M}_\sigma(-w;F) \,|dw| 
=\int_{\mathbb{C}} \widehat{\Phi}(w) \overline{\widetilde{M}_\sigma(w;F)} \,|dw| 
=\int_{\mathbb{C}} \Phi(z) M_\sigma(z;F) \,|dw| 
\end{gather*}
due to Parseval's identity, and therefore \eqref{eqM1} and \eqref{eqM1E} follow. 
The proof of the non-negativity of the function $M_\sigma(z;F)$ remains. 
For this, we assume $M_\sigma(z;F)<0$ for some region $U$. 
If we take $\Phi(z)$ as a non-negative function with a support included in $U$, then we have the contradiction. 
Due to the continuity of $M_\sigma(z;F)$, we see that $M_\sigma(z;F)$ is everywhere non-negative. 
\end{proof}

\section{Proof of Theorem \ref{thmM2}}\label{sec4}
We find that Theorem \ref{thmM1} imply Theorem \ref{thmM2} by the following lemma. 

\begin{lemma}\label{lemLes}
Let 
\begin{gather*}
K(x) 
=\left(\frac{\sin \pi x}{\pi x}\right)^2. 
\end{gather*}
Then for any $a,b \in \mathbb{R}$ with $a<b$, there exists a continuous function $F_{a,b}: \mathbb{R} \to \mathbb{R}$ such that the following conditions hold: for any $\omega>0$, 
\begin{enumerate}
\item[$\mathrm{(1)}$] 
$F_{a,b}(x)-1_{[a,b]}(x) \ll K(\omega(x-a))+K(\omega(x-b))$ for any $x \in \mathbb{R}$; 
\item[$\mathrm{(2)}$] 
$\displaystyle{ \int_{\mathbb{R}} (F_{a,b}(x)-1_{[a,b]}(x)) \,dx \ll \omega^{-1}}$; 
\item[$\mathrm{(3)}$] 
if $|x| \geq \omega$, then $\widehat{F}_{a,b}(x)=0$; 
\item[$\mathrm{(4)}$] 
$\widehat{F}_{a,b}(x) \ll (b-a)+\omega^{-1}$. 
\end{enumerate}
Here, 
\begin{gather*}
\widehat{F}_{a,b}(x) 
=\int_{\mathbb{R}} F_{a,b}(u) e^{ixu} \,|du| 
\end{gather*}
is the Fourier transformation of $F_{a,b}(x)$ with $|du|=(2\pi)^{-\frac{1}{2}} du$. 
\end{lemma}

\begin{proof}
This is Lemma 4.1 of \cite{Lester2014a} except for the difference of the definition of the Fourier transform, which does not affect the result. 
\end{proof}

\begin{proof}[Proof of Theorem \ref{thmM2}]
Again we consider only the case of $F \in \mathcal{S}_{\mathrm{I}}$. 
Assume that the rectangle $R$ is given as 
\begin{gather*}
R 
=\{ z=x+iy \in \mathbb{C} \mid a \leq x \leq b,~ c \leq y \leq d \}. 
\end{gather*}
Then we define for $z=x+iy \in \mathbb{C}$ 
\begin{gather}\label{eqPhi}
\Phi(z) 
=F_{a,b}(x) F_{c,d}(y). 
\end{gather}
We first find that the function $\Phi(z)$ belongs to the class $\Lambda$. 
The class $\Lambda$ is also written as 
\begin{gather*}
\Lambda
=\{ f \in L^1 \mid \text{$f$ is continuous and $\widehat{f} \in L^1$} \},
\end{gather*}
and hence we must check that $\Phi \in L^1$, $\Phi$ is continuous, and $\widehat{\Phi} \in L^1$. 
Since 
\begin{gather*}
\int_{\mathbb{C}} \Phi(z) \,|dz| 
=\int_{\mathbb{R}} F_{a,b}(x) \,|dx| \int_{\mathbb{R}} F_{c,d}(y) \,|dy|,
\end{gather*}
we see that $\Phi \in L^1$ by condition $\mathrm{(2)}$ of Lemma \ref{lemLes}. 
The function $\Phi(z)$ is continuous by its definition \eqref{eqPhi}, and furthermore, we have 
\begin{gather*}
\widehat{\Phi}(z) 
=\widehat{F}_{a,b}(x) \widehat{F}_{c,d}(y) 
=0 
\end{gather*}
if $|x|, |y| \geq \omega$ by condition $\mathrm{(3)}$. 
Thus also we have $\widehat{\Phi} \in L^1$. 
Therefore $\Phi(z)$ belongs to the class $\Lambda$, and we apply Theorem \ref{thmM1} for this function. 
Note that 
\begin{gather}\label{eqAppro}
\Phi(z)-1_R(z) 
\ll K(\omega(x-a))+K(\omega(x-b))+K(\omega(y-c))+K(\omega(y-d)) 
\end{gather}
by condition $\mathrm{(1)}$ of Lemma \ref{lemLes}. 
Then, let $\sigma>1-b^{-1}$ be fixed, and let $\theta, \delta>0$ with $\delta+3\theta>0$. 
We take $\omega=(\log T)^{\delta}$. 
Due to inequality \eqref{eqAppro}, Theorem \ref{thmM1} gives 
\begin{gather}\label{eqM2E}
\frac{1}{T} V_\sigma(T,R;F) 
=\int_{R} M_\sigma(z;F) \,|dz|
+E_1+E_2+E_3 
\end{gather}
for large $T$, where 
\begin{gather}\label{eqM2E1}
E_1 
\ll \exp\left(-\frac{1}{4} (\log T)^{\frac{2}{3}\theta} \right) \int_{\Omega} |\widehat{\Phi}(z)| \,|dz| 
+\int_{\mathbb{C} \setminus \Omega} |\widehat{\Phi}(z)| \,|dz|, 
\end{gather}
\begin{align}\label{eqM2E2}
E_2 
&\ll \frac{1}{T} \int_{0}^{T} K\left(\omega \left(\RE \frac{F'}{F}(\sigma+it) -a \right) \right) \,dt \\
&\qquad
+\frac{1}{T} \int_{0}^{T} K\left(\omega \left(\RE \frac{F'}{F}(\sigma+it)-b \right) \right) \,dt \nonumber \\
&\qquad\quad
+\frac{1}{T} \int_{0}^{T} K\left(\omega \left(\IM \frac{F'}{F}(\sigma+it)-c \right) \right) \,dt \nonumber \\
&\qquad\qquad
+\frac{1}{T} \int_{0}^{T} K\left(\omega \left(\IM \frac{F'}{F}(\sigma+it)-d \right) \right) \,dt, \nonumber 
\end{align}
and 
\begin{align}\label{eqM2E3}
E_3 
&\ll \int_{\mathbb{C}} K(\omega(x-a)) M_\sigma(z;F) \,|dz| 
+\int_{\mathbb{C}} K(\omega(x-b)) M_\sigma(z;F) \,|dz|  \\
&\qquad+\int_{\mathbb{C}} K(\omega(y-c)) M_\sigma(z;F) \,|dz| 
+\int_{\mathbb{C}} K(\omega(y-d)) M_\sigma(z;F)|dz|. \nonumber
\end{align}
All implied constants depend on $F, \sigma, \theta, \delta, \epsilon$. 
We estimate three error terms $E_1$, $E_2$, and $E_3$. 
The first term of the right hand side of \eqref{eqM2E1} is estimated as 
\begin{align*}
\exp\left(-\frac{1}{4} (\log T)^{\frac{2}{3}\theta} \right) \int_{\Omega} |\widehat{\Phi}(z)| \,|dz| 
&\ll\exp\left(-\frac{1}{4} (\log T)^{\frac{2}{3}\theta} \right) (\log T)^{2\delta} (b-a) (d-c) \\
&\ll (\log{T})^{-\delta} \nu_2(R) 
\end{align*}
for sufficiently large $T$ by condition $\mathrm{(4)}$ of Lemma \ref{lemLes}. 
We have 
\begin{gather*}
\int_{\mathbb{C} \setminus \Omega} |\widehat{\Phi}(z)| \,|dz| 
=0 
\end{gather*}
since $\widehat{\Phi}(z)=0$ if $|x|, |y| \geq \omega$. 
Therefore we obtain 
\begin{gather}\label{eqM2E1'}
E_1 
\ll \nu_2(R) (\log{T})^{-\delta}. 
\end{gather}
Next we estimate $E_2$. 
Since we have 
\begin{gather*}
K(\omega x) 
=\frac{2}{\omega^2} \int_{0}^{\omega} (\omega-u) \cos(2\pi xu) \,du 
=\frac{2}{\omega^2} \RE \int_{0}^{\omega} (\omega-u) e^{2\pi i xu} \,du, 
\end{gather*}
the first term of the right hand side of \eqref{eqM2E2} is estimated as 
\begin{align}\label{eqM2E2'}
&\frac{1}{T} \int_{0}^{T} K\left(\omega \left(\RE \frac{F'}{F}(\sigma+it)-a\right) \right) \,dt \\
&\ll \frac{1}{\omega^2} \int_{0}^{\omega} (\omega-u) 
\left| \frac{1}{T} \int_{0}^{T} \exp\left(2\pi iu \RE \frac{F'}{F}(\sigma+it) \right) \,dt \right| \,du. \nonumber
\end{align}
Proposition \ref{propM1} deduces 
\begin{gather*}
\frac{1}{T} \int_{0}^{T} \exp\left(2\pi iu \RE \frac{F'}{F}(\sigma+it) \right) \,dt
\ll \left|\widetilde{M}_\sigma(2\pi u;F)\right|
\end{gather*}
as $T\to\infty$, hence \eqref{eqM2E2'} is 
\begin{gather*}
\ll \frac{1}{\omega^2} \int_{0}^{\omega} (\omega-u) \left|\widetilde{M}_\sigma(2\pi u;F)\right| \,du 
\ll \frac{1}{\omega}
=(\log{T})^{-\delta}. 
\end{gather*}
The last inequality follows from Proposition \ref{propM3}. 
Since the reminder terms of \eqref{eqM2E2} are estimated in a similar way, we have  
\begin{gather}\label{eqM2E2''}
E_2 
\ll(\log{T})^{-\delta}. 
\end{gather}
The work of the estimate of $E_3$ remains. 
For this, we define 
\begin{gather*}
m_\sigma(x;F) 
=\int_{\mathbb{R}} M_\sigma(x+iy;F) \,|dy|.
\end{gather*}
Then the first term of the right hand side of \eqref{eqM2E3} is equal to 
\begin{gather*}
\int_{\mathbb{R}} K(\omega(x-a)) m_\sigma(x;F) \,|dx|.
\end{gather*}
The function $m_\sigma(x;F)$ is bounded on $\mathbb{R}$. 
In fact, it is continuous, and we see that 
\begin{gather*}
\int_{\mathbb{R}} m_\sigma(x;F) \,|dx| 
=\int_{\mathbb{C}} M_\sigma(x;F) \,|dz| 
=\widetilde{M}_\sigma(0;F) 
=1.
\end{gather*}
Therefore, we obtain 
\begin{gather*}
\int_{\mathbb{C}} K(\omega(x-a)) \mathcal{M}_\sigma(z;F) \,|dz| 
\ll \int_{\mathbb{R}} K(\omega(x-a)) \,dx 
\ll \frac{1}{\omega} 
=(\log{T})^{-\delta}. 
\end{gather*}
Estimating the remaining terms of \eqref{eqM2E3} similarly, we have 
\begin{gather}\label{eqM2E3'}
E_3 
\ll (\log{T})^{-\delta}. 
\end{gather}
By estimates \eqref{eqM2E1'}, \eqref{eqM2E2''}, and \eqref{eqM2E3'}, formula \eqref{eqM2E} gives 
\begin{gather*}
\frac{1}{T}V_\sigma(T,R;F)-\int_{R} M_\sigma(z;F) \,|dz| 
\ll (\nu_2(R)+1) (\log{T})^{-\delta}. 
\end{gather*}
Taking care of the assumption $\delta+3\theta<1/2$, we put $\theta=\epsilon/4$ and $\delta=1/2-\epsilon$ for arbitrarily small $\epsilon>0$. 
Then we obtain 
\begin{gather*}
(\nu_2(R)+1) (\log{T})^{-\delta} 
=(\nu_2(R)+1) (\log{T})^{-\frac{1}{2}+\epsilon}, 
\end{gather*}
which gives the result. 
\end{proof}


\providecommand{\bysame}{\leavevmode\hbox to3em{\hrulefill}\thinspace}
\providecommand{\MR}{\relax\ifhmode\unskip\space\fi MR }
\providecommand{\MRhref}[2]{%
  \href{http://www.ams.org/mathscinet-getitem?mr=#1}{#2}
}
\providecommand{\href}[2]{#2}

\end{document}